\documentclass[11pt]{amsart}
\usepackage{mathrsfs}
\usepackage{amssymb}

\usepackage{amscd}
\usepackage{amsmath, amssymb}
\usepackage{amsfonts}
\usepackage[colorlinks,linkcolor=blue,citecolor=blue, pdfstartview=FitH]{hyperref}
\usepackage[all]{xy}

\numberwithin{equation}{section}

\theoremstyle{plain}
\newtheorem{thm}{Theorem}[section]
\newtheorem{theorem}[thm]{Theorem}
\newtheorem{lemma}[thm]{Lemma}

\newtheorem{prop}[thm]{Proposition}

\theoremstyle{definition}

\newtheorem{remark}[thm]{Remark}

\newtheorem{definition}[thm]{Definition}

\newtheorem{example}[thm]{Example}

\newtheorem{defn-thm}[thm]{Definition-Theorem}

\newcommand{\btheorem}{\begin{theorem}}
\newcommand{\etheorem}{\end{theorem}}
\newcommand{\bproposition}{\begin{proposition}}
\newcommand{\eproposition}{\end{proposition}}
\newcommand{\bdefinition}{\begin{definition}}
\newcommand{\edefinition}{\end{definition}}
\newcommand{\bcorollary}{\begin{corollary}}
\newcommand{\ecorollary}{\end{corollary}}
\newcommand{\bproof}{\begin{proof}}
\newcommand{\eproof}{\end{proof}}

\newcommand{\eremark}{\end{remark}}
\newcommand{\eexample}{\end{example}}
\newcommand{\bexample}{\begin{example}}

\newcommand{\elemma}{\end{lemma}}
\newcommand{\blemma}{\begin{lemma}}

\renewcommand{\bar}{\overline}

\newcommand{\beq}{\begin{equation}}
\newcommand{\eeq}{\end{equation}}

\begin{document}
\title{Derivations on almost complex manifolds}
\author{Wei Xia}
\address{Wei Xia, Center of Mathematical Sciences, Zhejiang
University, Hangzhou, 310027, CHINA. } \email{xiaweiwei3@126.com}

\thanks{ }

\begin{abstract}
In this short note, we propose an unified method to derive formulas for derivations conjugated by exponential functions on an almost complex manifold.
\end{abstract}

\maketitle
%%%%%%%%%%%%%%%%%%%%%%%%%%%%%%%%%%%%%%%%%%%%%%%%%%%%%%%%%%%
%
%                 Section 1 Introduction
%
%%%%%%%%%%%%%%%%%%%%%%%%%%%%%%%%%%%%%%%%%%%%%%%%%%%%%%%%%%%
\section
{\bf Introduction}
%%%%%%%%%%%%%%%%%%%%%%%%%%%%%%%%%%%%%%%%%%%%%%%%%%%%%%%%%%%
Our main result of this note is the following
\begin{theorem}  Let $(M,J)$ be an almost complex manifold, $E$ be a smooth vector bundle on $M$, and $\nabla$ be a linear connection on $E$ with the decomposition according to bidegrees~:~$\nabla=\nabla^{1,0} + \nabla^{0,1} - i_{\theta} - i_{\bar{\theta}}$, where $\theta \in A^{2,0}(M,T^{0,1})$ is the torsion form of $J$. For any vector form $K\in A^{l+1}(M,TM),~l\geq 0$, denote the interior derivative by $i_K$ and the Lie derivatives by $\mathcal{L}_{K}:=[i_K,\nabla]$, $\mathcal{L}_{K}^{1,0}:=[i_K,\nabla^{1,0}]$, $\mathcal{L}_{K}^{0,1}:=[i_K,\nabla^{0,1}]$, respectively. For any $\phi\in A^{0,1}(M,T^{1,0})$ and $\psi\in A^{0,1}(M,T^{1,0})$, set  $e^{i_{\phi}}:=\sum_{k=0}^{\infty} \frac{i_{\phi}^k}{k!}$ and define $e^{i_{\bar{\phi}}}, e^{i_{\psi}}, e^{i_{\bar{\psi}}}$ in a similar way. Then we have
\begin{align*}
(1)~&e^{-i_{\phi}}\nabla e^{i_{\phi}}=\nabla - \mathcal{L}_{\phi} - \frac{1}{2}i_{[\phi,\phi]}- \frac{1}{3!}i_{[[\phi,\phi],\phi]^{\wedge}}~~;\\
(2)~&e^{-i_{\phi}}\nabla^{1,0} e^{i_{\phi}}= \nabla^{1,0}- \mathcal{L}_{\phi}^{1,0} - \frac{1}{2}i_{[\phi,\phi]^{A^{0,2}(T^{1,0})}}~~,\\
&e^{-i_{\phi}}\nabla^{0,1} e^{i_{\phi}}= \nabla^{0,1}- \mathcal{L}_{\phi}^{0,1} ~~;\\
(3)~&e^{-i_{\phi}}i_{\theta} e^{i_{\phi}}=i_{\theta}+i_{[\theta,\phi]^{\wedge}}+ \frac{1}{2}i_{[\theta,\phi]^{\wedge(2)}}+ \frac{1}{3!}i_{[\theta,\phi]^{\wedge(3)}}~~,\\
&e^{-i_{\phi}}i_{\bar{\theta}}e^{i_{\phi}}=i_{\bar{\theta}}~~;\\
(4)~&e^{-i_{\bar{\psi}}}i_{\phi}e^{i_{\bar{\psi}}}=i_{\phi+[\phi,\overline{\psi}]^{\wedge}+[\phi,\overline{\psi}]^{\wedge(2)}}~~,\\ &e^{-i_{\overline{\psi}}}i_{[\phi,\phi]}e^{i_{\overline{\psi}}}=i_{\sum_{j=0}^{j=3}\frac{1}{j!}[[\phi,\phi],\overline{\psi}]^{\wedge(j)}}~~;\\
(5)~&e^{-i_{\overline{\psi}}}\mathcal{L}_{\phi}e^{i_{\overline{\psi}}}=\mathcal{L}_{\phi-i_{\overline{\psi}}\phi}+
i_{\sum_{j=0}^{j=2}\frac{1}{(j+1)!}[[\phi,\overline{\psi}],\overline{\psi}]^{\wedge(j)}}-
i_{\sum_{j=0}^{j=2}\frac{1}{(j+2)!}[[i_{\overline{\psi}}\phi,\overline{\psi}],\overline{\psi}]^{\wedge(2)}}~~;\\
(6)~&e^{-i_{\overline{\psi}}}e^{-i_{\phi}}\nabla e^{i_{\phi}}e^{i_{\overline{\psi}}}
    =\nabla-\mathcal{L}_{\overline{\psi}}-\frac{1}{2}i_{[\overline{\psi},\overline{\psi}]}- \frac{1}{3!}i_{[[\bar{\psi},\bar{\psi}],\bar{\psi}]^{\wedge}} -\mathcal{L}_{\phi-i_{\overline{\psi}}\phi}-
i_{\sum_{j=0}^{j=2}\frac{1}{(j+1)!}[[\phi,\overline{\psi}],\overline{\psi}]^{\wedge(j)}} \\
&+i_{\sum_{j=0}^{j=2}\frac{1}{(j+2)!}[[i_{\overline{\psi}}\phi,\overline{\psi}],\overline{\psi}]^{\wedge(2)}}-
\frac{1}{2}i_{\sum_{j=0}^{j=3}\frac{1}{j!}[[\phi,\phi],\overline{\psi}]^{\wedge(j)}} -\frac{1}{3!}i_{\sum_{j=0}^{j=3}\frac{1}{j!}[[[\phi,\phi],\phi]^{\wedge},\bar{\psi}]^{\wedge(j)}}~~.
\end{align*}
\end{theorem}

$\mathbf{NOTE:}$~The corresponding formulas in a previous version of this paper~\cite{Xia19dp} are not correct. We mistakenly take $[A^{0,1}(M,T^{1,0}),A^{0,1}(M,T^{1,0})]\in A^{0,2}(M,T^{1,0})$ as granted which does not hold for general almost complex manifolds, see \eqref{notice}. The correct form of formulas $(1), (2), (3)$ are first obtained by Fu-Liu in~\cite{FL19}. The main contribution of this note is thus to provide an alternative proof of these formulas.

If $J$ is integrable, then $[A^{0,1}(M,T^{1,0}),A^{0,1}(M,T^{1,0})]\in A^{0,2}(M,T^{1,0})$ holds and $(1), (2)$ is reduced to the extension formulas proved in~\cite{LRY15}. After reviewing some basic facts about derivations in Section \ref{sec rieview}, we prove our main result in Section \ref{Derivations on an ac manifold}.

%
%%%%%%%%%%%%%%%%%%%%%%%%%%%%%%%%%%%%%%%%%%%%%%%%%%%%%%%%%%%
\section{Derivations and bracket operations on a real manifold }
\label{sec rieview}
In this section we review some basic facts about derivations and refer the readers to chapter II of \cite{KMS93} for more details.
Let $M$ be a smooth manifold of dimension $n$, $A(M)=\bigoplus_{k=0}^{n} A^k(M)$ be its exterior algebra of differential forms. A (graded) \emph{derivation} $D$ on $A(M)$ is a $\mathbb{R}$-linear map $D:A(M)\to A(M)$ with $D(A^l(M))\subseteq A^{l+k}(M)$ and $D(\xi\wedge\eta)=(D\xi)\wedge\eta+(-1)^{kl}\xi\wedge(D\eta)$ for $\xi\in A^l(M)$.The integer $k$ is called the degree of $D$. We denote by $D^k(M)$ the space of all derivations of degree $k$ on $A(M)$. For any $D_1\in D^{k_1}(M)$ and $D_2\in D^{k_2}(M)$, the graded commutator is defined by $[D_1,D_2]:=D_1D_2-(-1)^{k_1k_2}D_2D_1.$ With respect to this bracket operation, the space of all derivations $D(M)=\bigoplus_kD^k(M)$ becomes a graded Lie algebra.
\subsection{The interior derivative and Lie derivative}
For a smooth vector field $X$ on $M$, we have the interior derivative $i_X$ which is a derivation of degree $-1$. For a vector valued $(k+1)$-form $K\in A^{k+1}(M,TM)$, we can associate a derivation $i_K$ of degree $k$ by setting $i_K\varphi:=\xi\wedge (i_X\varphi)$ , if $K=\xi\otimes X$ for a $(k+1)$-form $\xi$ and a vector field $X$, where $\varphi\in A(M)$. The exterior derivative $d$ is a derivation of degree $1$. The Lie derivative $\mathcal{L}_X:=[i_X,d]=i_Xd+di_X$ is a derivation of degree $0$. Similarly, the Lie derivative $\mathcal{L}_K:=[i_K,d]=i_Kd-(-1)^{k}di_K$, where $K\in A^{k+1}(M,TM)$, is a derivation of degree $k+1$. In fact, for any $D\in D^k(M)$, there exist unique $K\in A^{k}(M,TM)$ and $L\in A^{k+1}(M,TM)$ such that \begin{equation}\label{D=lK+iL}
D=\mathcal{L}_K+i_L,
\end{equation}
and $L=0$ if and only if $[D,d]=0$, $K=0$ if and only if $D$ is algebraic.
\subsection{Algebraic derivation and Nijenhuis-Richardson bracket}
A derivation $D$ on $A(M)$ is called \emph{algebraic} if $Df=0, \forall f\in A^0(M)$. Every algebraic derivation of degree $k$ on $A(M)$ is of the form $i_K$ for some unique $K\in A^{k+1}(M,TM)$. For any two algebraic derivation $i_K\in D^k(M)$ and $i_L\in D^l(M)$, $[i_K,i_L]$ is again an algebraic derivation , hence $[i_K,i_L]=i_{[K,L]^{\wedge}}$ for some unique $[K,L]^{\wedge}\in A^{k+l+1}(M,TM)$. We have
\begin{equation}\label{ii commutator}
[i_K,i_L]=i_{[K,L]^{\wedge}}=i_{({i_KL-(-1)^{kl}i_LK})}.
\end{equation}
The operation $[\cdot,\cdot]^{\wedge}$ is called the \emph{Nijenhuis-Richardson bracket}. $i_KL$ is defined as $i_KL:=\xi\wedge(i_X \eta)\otimes Y$ for $K=\xi\otimes X$ and $L=\eta\otimes Y$.

\subsection{The Fr\"{o}licher-Nijenhuis bracket}

For any $K\in A^{k}(M,TM)$, $L\in A^{l}(M,TM)$, $[\mathcal{L}_K,\mathcal{L}_L]$ is a derivation of degree $k+l$ such that $[[\mathcal{L}_K,\mathcal{L}_L],d]=0$, hence $[\mathcal{L}_K,\mathcal{L}_L]=\mathcal{L}_{[K,L]}$ for some unique $[K,L]\in A^{k+l}(M,TM)$. This operation $[\cdot,\cdot]$ is called the \emph{Fr\"{o}licher-Nijenhuis bracket}. For $k=l=0$, this bracket coincides with the Lie bracket of vector fields. On an almost complex manifold $(M,J)$, the Newlander-Nirenberg theorem says that $J$ is integrable if and only if $[J,J]=0$, where $J:TM\to TM$ is considered as a vector $1$-form. On a complex manifold, the Fr\"{o}licher-Nijenhuis bracket can be extended $\mathbb{C}$-linearly and is exactly the bracket operation appeared in Kodaira-Spencer's deformation theory \cite{MK71}. Nevertheless, we should notice a vital difference between the integrable case and the general case: for $\phi\in A^{0,1}(M,T^{1,0})$ and $\psi\in A^{0,1}(M,T^{1,0})$, $[\phi,\psi]\in A^{0,2}(M,T^{1,0})$ does not always hold if $J$ is not integrable. In fact, we have\footnote{I owed this observation to Hai-Sheng Liu, see Lemma 3.3 of~\cite{FL19} and Theorem 8.7, (6) of~\cite[pp.\,70]{KMS93}.}
\begin{equation}\label{notice}
[\phi,\psi]\in A^{0,2}(M,T^{1,0})\oplus A^{1,1}(M,T^{1,0})\oplus A^{0,2}(M,T^{0,1}).
\end{equation}
We denote by $[\phi,\psi]^{A^{0,2}(T^{1,0})}$ the component of $[\phi,\psi]$ which lies in $A^{0,2}(M,T^{1,0})$. Similarly, we write $[\phi,\psi]^{A^{1,1}(T^{1,0})}$ and $[\phi,\psi]^{A^{0,2}(T^{0,1})}$ for the other components.
\subsection{A useful commutator relation}

The following commutator relation will be useful for our purpose: for $K\in A^{k}(M,TM)$, $L\in A^{l+1}(M,TM)$, we have
\begin{equation}\label{li commutator}
[\mathcal{L}_K,i_{L}]=i_{[K,L]}-(-1)^{kl}\mathcal{L}_{i_LK}.
\end{equation}
See \cite{Mic86}, \cite{LR11} and \cite{LRY15} for various forms and generalizations of this formula.
%%%%%%%%%%%%%%%%%%%%%%%%%%%%%%%%%%%%%%%%%%%%%%%%%%%%%%%%%%%
%
%   Section 3 Current equation on the compactification of the moduli space of polarized Calabi-Yau manifolds
%
%%%%%%%%%%%%%%%%%%%%%%%%%%%%%%%%%%%%%%%%%%%%%%%%%%%%%%%%%%%
\section{Derivations on an almost complex manifold }
\label{Derivations on an ac manifold}
Now, let $(M,J)$ be an almost complex manifold of real dimension $2n$, its complexified tangent bundle $T_{\mathbb{C}}(M)=T_{\mathbb{R}}(M)\otimes \mathbb{C}$ admits a decomposition $T_{\mathbb{C}}(M)=T^{1,0}\oplus T^{0,1}$. For each $k\leq 0$, the space of (complexified) $k$-forms $A_{\mathbb{C}}^k(M):=A^k(M)\otimes \mathbb{C}$ admits a decomposition $A^k(M)=\bigoplus_{p+q=k} A^{p,q}(M)$, where $A^{p,q}(M)=\wedge^pT^{*1,0}\otimes\wedge^qT^{*0,1}$. Hence $A_{\mathbb{C}}(M):=A(M)\otimes \mathbb{C}=A(E)$. In what follows, we omit the subscript "$\mathbb{C}$" and make the convention that differential forms are always complex valued unless otherwise stated. Replacing $\mathbb{R}$-linearity by $\mathbb{C}$-linearity, the notion of graded derivation on $A(M)$ is similarly defined. Since now $A(M)$ has a bigrading structure, we can make a refinement. A \emph{bigraded derivation} of bidegree $(k,l)$ on $A(M)$ is a $\mathbb{C}$-linear map $D:A(M)\to A(M)$ with $D(A^{p,q}(M))\subseteq A^{p+k,q+l}(M)$ and $D(\xi\wedge\eta)=(D\xi)\wedge\eta+(-1)^{(k+l)m}\xi\wedge(D\eta)$ for $\xi\in A^m(M)$. By definition, A bigraded derivation of bidegree $(k,l)$ is necessarily a graded derivation of degree $k+l$. The space of all bigraded derivation of bidegree $(k,l)$ on $A(M)$ is denoted by $D^{k,l}(M)$, then $D^k(M)=\bigoplus_{p+q=k} D^{p,q}(M)$. Note that, by extending $\mathbb{C}$-linearly, all the constructions in previous section can be applied in the present situation.
\begin{example}[c.f. chapter VIII of \cite{Dem12}]\label{decomposition of d}
 The exterior derivative $d$ admit a decomposition into $4$ bigraded derivations:
\begin{equation*}d=\partial + \bar{\partial} - i_{\theta} - i_{\bar{\theta}},
\end{equation*}
where $\partial:=\sum_{p,q}\Pi^{p+1,q}d\Pi^{p,q}$ with $\Pi^{p,q}$ being the projection $A(M)\to A^{p,q}(M)$ and $\bar{\partial}:=\sum_{p,q}\Pi^{p,q+1}d\Pi^{p,q}$. And $\theta \in A^{2,0}(M,T^{0,1})$ is the torsion form of $J$ which is defined by $\theta(X,Y):=[X,Y]^{0,1}$ for $X,Y \in A^{0}(M,T^{1,0})$, where $[X,Y]^{0,1}$ is the $(0,1)$ part of the vector field $[X,Y]$. We see that the bidegree of $\partial,\bar{\partial},i_{\theta}, i_{\bar{\theta}}$ are $(1,0),(0,1),(2,-1),(-1,2)$ respectively.
\end{example}

\begin{example}
For $\phi\in A^{0,k}(M,T^{1,0})$, define
\begin{equation*}\mathcal{L}_{\phi}^{1,0}:=[i_\phi,\partial]=i_\phi\partial-\partial i_\phi~~~~~~~~ \text{and}~~~~~~~~ \mathcal{L}_{\phi}^{0,1}:=[i_\phi,\overline{\partial}]=i_\phi\overline{\partial}-\overline{\partial} i_\phi,
\end{equation*}
then  $\mathcal{L}_{\phi}^{1,0}\in D^{0,k}(M)$ and $\mathcal{L}_{\phi}^{0,1}\in D^{-1,k+1}(M)$.
\end{example}

It is clear that if $D\in D^{p,q}(M)$ is algebraic, $D=i_{L}$ for some unique $L\in A^{p+1,q}(M,T^{1,0})\oplus A^{p,q+1}(M,T^{0,1})$. We can also make the following refinement of (\ref{D=lK+iL}):
\begin{prop}\label{D=l10K+l01K+iL}
Let $D\in D^{p,q}(M)$, then we have
\begin{equation}
D=\mathcal{L}^{1,0}_{K^{'p,q}}+\mathcal{L}^{0,1}_{K^{''p,q}}+i_{L^{'p+1,q}}+i_{L^{''p,q+1}},
\end{equation}
for some $K^{'p,q}\in A^{p,q}(M,T^{1,0})$, $K^{''p,q}\in A^{p,q}(M,T^{0,1})$ and $L^{'p,q}\in A^{p,q}(M,T^{1,0})$, $L^{''p,q}\in A^{p,q}(M,T^{0,1})$.
\end{prop}
\begin{proof}Let $D\in D^{p,q}(M)$, then $D=\mathcal{L}_K+i_L,$ for some unique $K\in A^{p+q}(M,TM)$ and $L\in A^{p+q+1}(M,TM)$. We can write
\begin{equation*}
D=\sum_{a+b=p+q}(\mathcal{L}^{1,0}_{K^{'a,b}}+\mathcal{L}^{1,0}_{K^{''a,b}}+\mathcal{L}^{0,1}_{K^{'a,b}}+\mathcal{L}^{0,1}_{K^{''a,b}})
+\sum_{a+b=p+q+1}i_{L^{'a,b}+L^{''a,b}},
\end{equation*}
 where $\forall a,b, K^{'a,b}\in A^{a,b}(M,T^{1,0})$, $K^{''a,b}\in A^{a,b}(M,T^{0,1})$ and similarly $L^{'a,b}\in A^{a,b}(M,T^{1,0})$, $L^{''a,b}\in A^{a,b}(M,T^{0,1})$. It follows from $D^k(M)=\bigoplus_{p+q=k} D^{p,q}(M)$ that
\begin{equation*}
D=\mathcal{L}^{1,0}_{K^{'p,q}}+\mathcal{L}^{1,0}_{K^{''p-1,q+1}}+\mathcal{L}^{0,1}_{K^{'p+1,q-1}}+\mathcal{L}^{0,1}_{K^{''p,q}}+i_{L^{'p+1,q}}+i_{L^{''p,q+1}}.
\end{equation*}
It is clear that $\mathcal{L}^{1,0}_{K^{''p-1,q+1}}$ and $\mathcal{L}^{0,1}_{K^{'p+1,q-1}}$ are algebraic so that $\mathcal{L}^{1,0}_{K^{''p-1,q+1}}=i_{R^{''p,q+1}}$, $\mathcal{L}^{0,1}_{K^{'p+1,q-1}}=i_{R^{'p+1,q}}$ for some $R^{''p,q+1}\in A^{p,q+1}(M,T^{0,1})$ and $R^{'p+1,q}\in A^{p+1,q}(M,T^{1,0})$.
\end{proof}
\begin{remark}It is important to notice that the uniqueness part of (\ref{D=lK+iL}) is lost.
\end{remark}

%We want to derive explicit formula for the operator \bar{\partial}.
 Let $\phi\in A^{0,1}(M,T^{1,0})$, then $i_{\phi}$ is nilpotent : $(i_{\phi})^{n+1}\xi=0, \forall \xi \in A(M)$, so that the operator
\begin{equation*} e^{i_{\phi}}:=\sum_{k=0}^{\infty} \frac{i_{\phi}^k}{k!}~~~~~~: A(M) \longrightarrow A(M)
\end{equation*}
is well-defined. Since $e^{i_{\phi}}e^{-i_{\phi}}=e^{-i_{\phi}}e^{i_{\phi}}=e^0$ is the identity operator, $e^{-i_{\phi}}$ is the inverse operator of $e^{i_{\phi}}$.

\begin{definition}Let $R$ be an unitary associative algebra over $\mathbb{Q}$ (not necessarily commutative). For any $x,y\in R$, we say that $x$ is \emph{finitely commutable} with $y$ if there is a positive integer $k$ such that
\begin{equation}\label{finitly commu}\underbrace{[\cdots[}_{k~~\text{times}}x,\overbrace{y],y],\cdots,y]}^{k~~\text{times}}=0,
\end{equation}
where $[x,y]=xy-yx$ is the usual commutator. If $x$ is finitely commutable with $y$, the least integer $k$ such that (\ref{finitly commu}) holds is called the \emph{commutable degree} of $(x,y)$, and in this case we say $x$ is \emph{$k$-commutable} with $y$. We will simply denote the $k$ times bracket in (\ref{finitly commu}) by $[x,y]^{(k)}$ and make the convention that $[x,y]^{(0)}:=x$.
\end{definition}

\subsection{A commutator lemma}
The following lemma is perhaps well-known to experts, see \cite[Lem.\,2.7]{Got05} and \cite[pp.\,66,~Exercise~V.1.]{Ma05}. For the readers' convenience, we will present here two different proofs.
\begin{lemma}\label{commutator lemma}
Let $R$ be an unitary associative algebra over $\mathbb{Q}$ (not necessarily commutative), $y\in R$ be a nilpotent element. Assume $x\in R$ is $k$-commutable with $y$, i.e. $[x,y]^{(k)}=0$, then
\begin{equation}\label{e commutator} e^{-y}xe^{y}=\sum_{i=0}^{i=k-1}\frac{1}{i!}[x,y]^{(i)},
\end{equation}
where $e^{y}:=1+y+\frac{y^2}{2!}+\cdots$ is the exponential function.
\end{lemma}

\begin{proof}First, we set
\[f(t):=e^{-ty}xe^{ty},
\]
where $t$ is a real variable. It can be proved inductively that
\[\frac{d^k f}{dt^k}(0)=[x,y]^{(k)}.
\]
Hence, we have
\[e^{-y}xe^{y}=f(1)=\sum_k\frac{d^k f}{dt^k}(0)\frac{1}{k!}=\sum_k\frac{1}{k!}[x,y]^{(k)}.
\]
\end{proof}

\begin{proof}Assume $y^l=0$ for some positive integer $l$ and set $N=\max\{k,l\}$, then
\begin{align*}xe^{y}&=x(1+y+\frac{y^2}{2!}+\cdots)\\
                                 &=x+xy(1+\frac{y}{2!}+\frac{y^2}{3!}+\cdots)\\
                                 &=x+([x,y]+yx)(1+\frac{y}{2!}+\frac{y^2}{3!}+\cdots)\\
                                 &=(1+y)x+[x,y]+([x,y]+yx)(\frac{y}{2!}+\frac{y^2}{3!}+\cdots)\\
                                 &=(1+y)x+[x,y]+([x,y]^{(2)}+y[x,y]+y[x,y]+y^2x)(\frac{1}{2!}+\frac{y}{3!}+\cdots)\\
                                 &=\sum_{i=0}^{i=2}\frac{y^i}{i!}x+(1+y)[x,y]+\frac{1}{2}[x,y]^{(2)}+([x,y]^{(2)}+2y[x,y]+y^2x)(\frac{y}{3!}+\cdots)\\
                                 &=\sum_{i=0}^{i=3}\frac{y^i}{i!}x+\sum_{i=0}^{i=2}\frac{y^i}{i!}[x,y]+\frac{1}{2}(1+y)[x,y]^{(2)}+\frac{1}{3!}[x,y]^{(3)}\\
                                 & +([x,y]^{(3)}+3y[x,y]^{(2)}+3y^2[x,y]+y^3x)(\frac{y}{4!}+\cdots)\\
                                 &=\sum_{i=0}^{i=4}\frac{y^i}{i!}x+\sum_{i=0}^{i=3}\frac{y^i}{i!}[x,y]+\frac{1}{2}\sum_{i=0}^{i=2}\frac{y^i}{i!}[x,y]^{(2)}+\frac{1}{3!}(1+y)[x,y]^{(3)}+\frac{1}{4!}[x,y]^{(4)}\\
                                 & +([x,y]^{(4)}+4y[x,y]^{(3)}+6y^2[x,y]^{(2)}+4y^3[x,y]+y^4x)(\frac{y}{5!}+\cdots)\\
                                 &=\cdots=\sum_{i=0}^{i=2N}\frac{y^i}{i!}x+\sum_{i=0}^{i=2N-1}\frac{y^i}{i!}[x,y]+\cdots+\frac{1}{N!}\sum_{i=0}^{i=N}\frac{y^i}{i!}[x,y]^{(N)}+\underbrace{\cdots+\frac{1}{(2N)!}[x,y]^{(2N)}}_{=0}\\
                                 & +\sum_{i=0}^{i=2N} {2N \choose i}y^i[x,y]^{(2N-i)}(\frac{y}{(2N+1)!}+\cdots),\\
\end{align*}
where ${2N \choose i}$ are the binomial coefficients. Now, by our assumption, we know that
\begin{equation*}\sum_{i=0}^{i=2N} {2N \choose i} y^i[x,y]^{(2N-i)}=0~~\text{and}~~
\sum_{i=0}^{i=2N}\frac{y^i}{i!}=\sum_{i=0}^{i=2N-1}\frac{y^i}{i!}=\cdots=\sum_{i=0}^{i=N}\frac{y^i}{i!}=e^y,
\end{equation*}
thus $xe^{y}=e^{y}\sum_{i=0}^{i=k-1}\frac{1}{i!}[x,y]^{(i)}~\Rightarrow~e^{-y}xe^{y}=\sum_{i=0}^{i=k-1}\frac{1}{i!}[x,y]^{(i)}$.
\end{proof}

\subsection{Derivations on the algebra of vector bundle valued forms}
Let $E$ be a smooth vector bundle on the the almost complex manifold $(M,J)$ and $\nabla$ be a linear connection on $E$. The space of E-valued differential forms $A(E)$ can be decomposed as $A(E)=\bigoplus_{p,q}A^{p,q}(E)$. From the work of Michor~\cite{Mic86}, we know that a similar theory as those described in Section \ref{sec rieview} holds in this setting, in particular, formulas (\ref{ii commutator}) and (\ref{li commutator}) are valid with the Lie derivative defined by $\mathcal{L}_K:=[i_K,\nabla]=i_K\nabla-(-1)^{k}\nabla i_K$, where $K\in A^{k+1}(M,TM)$. See \cite[Th.\,3.16]{Mic86} for the proof of (\ref{li commutator}). The space of all bigraded derivation of bidegree $(k,l)$ on $A(E)$ is denoted by $D^{k,l}(E)$, then $D^k(E)=\bigoplus_{p+q=k} D^{p,q}(E)$.

As in Example \ref{decomposition of d}, the connection $\nabla$ admits a decomposition\footnote{This can be checked easily. In fact, for $u\in A^{0}(M,E)$ this is clear; for $u\in A^{1}(M,E)$, write locally $u=\alpha_i\otimes s_i$, where $\alpha_i$ are $1$-forms and $s_i$ is a local smooth frame, then we have $\nabla u=\nabla(\alpha_i\otimes s_i)=d\alpha_i\otimes s_i- \alpha_i\otimes (\nabla^{1,0} + \nabla^{0,1})s_i =(\nabla^{1,0} + \nabla^{0,1} - i_{\theta} - i_{\bar{\theta}})u$.}:
\begin{equation*}
\nabla=\nabla^{1,0} + \nabla^{0,1} - i_{\theta} - i_{\bar{\theta}},
\end{equation*}
where $\nabla^{1,0}:=\sum_{p,q}\Pi^{p+1,q}\nabla\Pi^{p,q}$ with $\Pi^{p,q}$ being the projection $A(E)\to A^{p,q}(E)$ and $\nabla^{0,1}:=\sum_{p,q}\Pi^{p,q+1}\nabla\Pi^{p,q}$. Define $\mathcal{L}_{K}^{1,0}:=[i_K,\nabla^{1,0}]=i_K\nabla^{1,0}-\nabla^{1,0} i_K$ and $\mathcal{L}_{K}^{0,1}:=[i_K,\nabla^{0,1}]=i_K\nabla^{0,1}-\nabla^{0,1} i_K$ as usual. We make two observations. For $\phi\in A^{0,k}(M,T^{1,0})$, $\psi\in A^{0,l}(M,T^{1,0})$, it follows easily form the definition that
\begin{equation}\label{ii same type commutator} i_{\phi}i_{\psi}=(-1)^{(k+1)(l+1)}i_{\psi}i_{\phi},
\end{equation}
and since $i_{\psi}\phi=0$, by (\ref{li commutator}) we have
\begin{equation}\label{li same type commutator}
[\mathcal{L}_{\phi},i_{\psi}]=i_{[\phi,\psi]}.
\end{equation}

\begin{lemma}\label{lie1,0-i commutator} For $\phi\in A^{0,1}(M,T^{1,0})$, $\psi\in A^{0,1}(M,T^{1,0})$, we have
\begin{itemize}
  \item[1.]
$[\mathcal{L}^{1,0}_{\phi},i_{\psi}]=i_{[\phi,\psi]^{A^{0,2}(T^{1,0})}}$~;
  \item[2.]
$[\mathcal{L}^{0,1}_{\phi},i_{\psi}]=0$~;
  \item[3.]
$-[[\phi,\theta]^\wedge, \psi]^\wedge=[\phi,\psi]^{A^{1,1}(T^{1,0})}+[\phi,\psi]^{A^{0,2}(T^{0,1})}$~.
\end{itemize}
\end{lemma}

\begin{proof}Since
\[
\mathcal{L}_{\phi}= \mathcal{L}_{\phi}^{1,0}+ \mathcal{L}^{0,1}_{\phi}- i_{[\phi,\theta]^\wedge}\in D^{0,1}(E)\oplus D^{-1,2}(E)\oplus D^{1,0}(E),
\]
where $[\phi,\theta]^\wedge\in A^{1,1}(M,T^{0,1})\oplus A^{2,0}(M,T^{1,0})$, we have
\[
[\mathcal{L}_{\phi}, i_{\psi}]= [\mathcal{L}^{1,0}, i_{\psi}]+ [\mathcal{L}^{0,1}_{\phi}, i_{\psi}]- [i_{[\phi,\theta]^\wedge}, i_{\psi}]\in D^{-1,2}(E)\oplus D^{-2,3}(E)\oplus D^{0,1}(E)~.
\]
The conclusion then follows from \eqref{li same type commutator} and \eqref{notice}.
\end{proof}

As in (\ref{finitly commu}), we use $[x,y]^{\wedge(k)}$ to denote the $k$-times Nijenhuis-Richardson bracket of $x$ with $y$, and $[x,y]^{\wedge(0)}:=x.$

\begin{theorem} Let $\phi\in A^{0,1}(M,T^{1,0})$ and $\psi\in A^{0,1}(M,T^{1,0})$, then we have

\begin{align*}
(1)~&e^{-i_{\phi}}\nabla e^{i_{\phi}}=\nabla - \mathcal{L}_{\phi} - \frac{1}{2}i_{[\phi,\phi]}- \frac{1}{3!}i_{[[\phi,\phi],\phi]^{\wedge}}~~;\\
(2)~&e^{-i_{\phi}}\nabla^{1,0} e^{i_{\phi}}= \nabla^{1,0}- \mathcal{L}_{\phi}^{1,0} - \frac{1}{2}i_{[\phi,\phi]^{A^{0,2}(T^{1,0})}}~~,\\
&e^{-i_{\phi}}\nabla^{0,1} e^{i_{\phi}}= \nabla^{0,1}- \mathcal{L}_{\phi}^{0,1} ~~;\\
(3)~&e^{-i_{\phi}}i_{\theta} e^{i_{\phi}}=i_{\theta}+i_{[\theta,\phi]^{\wedge}}+ \frac{1}{2}i_{[\theta,\phi]^{\wedge(2)}}+ \frac{1}{3!}i_{[\theta,\phi]^{\wedge(3)}}~~,\\
&e^{-i_{\phi}}i_{\bar{\theta}}e^{i_{\phi}}=i_{\bar{\theta}}~~;\\
(4)~&e^{-i_{\bar{\psi}}}i_{\phi}e^{i_{\bar{\psi}}}=i_{\phi+[\phi,\overline{\psi}]^{\wedge}+[\phi,\overline{\psi}]^{\wedge(2)}}~~,\\ &e^{-i_{\overline{\psi}}}i_{[\phi,\phi]}e^{i_{\overline{\psi}}}=i_{\sum_{j=0}^{j=3}\frac{1}{j!}[[\phi,\phi],\overline{\psi}]^{\wedge(j)}}~~;\\
(5)~&e^{-i_{\overline{\psi}}}\mathcal{L}_{\phi}e^{i_{\overline{\psi}}}=\mathcal{L}_{\phi-i_{\overline{\psi}}\phi}+
i_{\sum_{j=0}^{j=2}\frac{1}{(j+1)!}[[\phi,\overline{\psi}],\overline{\psi}]^{\wedge(j)}}-
i_{\sum_{j=0}^{j=2}\frac{1}{(j+2)!}[[i_{\overline{\psi}}\phi,\overline{\psi}],\overline{\psi}]^{\wedge(2)}}~~;\\
(6)~&e^{-i_{\overline{\psi}}}e^{-i_{\phi}}\nabla e^{i_{\phi}}e^{i_{\overline{\psi}}}
    =\nabla-\mathcal{L}_{\overline{\psi}}-\frac{1}{2}i_{[\overline{\psi},\overline{\psi}]}- \frac{1}{3!}i_{[[\bar{\psi},\bar{\psi}],\bar{\psi}]^{\wedge}} -\mathcal{L}_{\phi-i_{\overline{\psi}}\phi}-
i_{\sum_{j=0}^{j=2}\frac{1}{(j+1)!}[[\phi,\overline{\psi}],\overline{\psi}]^{\wedge(j)}} \\
&+i_{\sum_{j=0}^{j=2}\frac{1}{(j+2)!}[[i_{\overline{\psi}}\phi,\overline{\psi}],\overline{\psi}]^{\wedge(2)}}-
\frac{1}{2}i_{\sum_{j=0}^{j=3}\frac{1}{j!}[[\phi,\phi],\overline{\psi}]^{\wedge(j)}} -\frac{1}{3!}i_{\sum_{j=0}^{j=3}\frac{1}{j!}[[[\phi,\phi],\phi]^{\wedge},\bar{\psi}]^{\wedge(j)}}~~.
\end{align*}

\end{theorem}

\begin{proof}All these follows easily from Lemma \ref{commutator lemma}.

For (1), note that $[\nabla,i_{\phi}]=-\mathcal{L}_{\phi}$. By (\ref{ii same type commutator}) and (\ref{li same type commutator}),
\begin{equation*}
[\nabla,i_{\phi}]^{(2)}=[-\mathcal{L}_{\phi},i_{\phi}]=-i_{[\phi,\phi]}~~~ \text{and}~~~[\nabla,i_{\phi}]^{(3)}=-[i_{[\phi,\phi]},i_{\phi}]=-i_{[[\phi,\phi],\phi]^{\wedge}}.
\end{equation*}
From \eqref{ii commutator} and \eqref{notice} we see that
\[
[[\phi,\phi],\phi]^{\wedge}\in A^{0,2}(M,T^{1,0})~~~ \text{and}~~~\big[[[\phi,\phi],\phi]^{\wedge},\phi\big]^{\wedge}=0,
\]
so
\[
[\nabla,i_{\phi}]^{(4)}=-i_{[[\phi,\phi],\phi]^{\wedge(2)}}=-i_{\big[[[\phi,\phi],\phi]^{\wedge},\phi\big]^{\wedge}}=0.
\]

For (2), by Lemma \ref{lie1,0-i commutator}, we have
\begin{equation*}
[\nabla^{1,0},i_{\phi}]=-\mathcal{L}_{\phi}^{1,0},~[\nabla^{1,0},i_{\phi}]^{(2)}=[-\mathcal{L}_{\phi}^{1,0},i_{\phi}]=-i_{[\phi,\phi]^{A^{0,2}(T^{1,0})}},~[\nabla^{1,0},i_{\phi}]^{(3)}=0
\end{equation*}
and
\begin{equation*}
[ \nabla^{0,1},i_{\phi}]=-\mathcal{L}_{\phi}^{0,1},~[ \nabla^{0,1},i_{\phi}]^{(2)}=[-\mathcal{L}_{\phi}^{0,1},i_{\phi}]=0.
\end{equation*}
For (3), we note that
\begin{align*}
[i_{\theta},i_{\phi}]&=i_{[\theta,\phi]^{\wedge}},~\text{where}~[\theta,\phi]^{\wedge}\in A^{1,1}(M,T^{0,1})\oplus A^{2,0}(M,T^{1,0}),\\
[i_{\theta},i_{\phi}]^{(2)}&=i_{[\theta,\phi]^{\wedge(2)}},~\text{where}~[\theta,\phi]^{\wedge(2)}\in A^{0,2}(M,T^{0,1})\oplus A^{1,1}(M,T^{1,0}),\\
[i_{\theta},i_{\phi}]^{(3)}&=i_{[\theta,\phi]^{\wedge(3)}},~\text{where}~[\theta,\phi]^{\wedge(3)}\in A^{0,2}(M,T^{1,0}),
\end{align*}
and so that $[i_{\theta},i_{\phi}]^{(4)}=i_{[\theta,\phi]^{\wedge(4)}}=0$. Also, $[i_{\bar{\theta}},i_{\phi}]=i_{[\bar{\theta},\phi]^{\wedge}}=0$.

For (4),  we note that $[i_{\phi},i_{\overline{\psi}}]^{(2)}=i_{[\phi,\overline{\psi}]^{\wedge(2)}}$ and $[\phi,\overline{\psi}]^{\wedge(2)}\in A^{1,0}(M,T^{0,1})$ implies that $[i_{\phi},i_{\overline{\psi}}]^{(3)}=[[i_{\phi},i_{\overline{\psi}}]^{(2)},i_{\overline{\psi}}]=[i_{[\phi,\overline{\psi}]^{\wedge(2)}},i_{\overline{\psi}}]=0$. Similarly,
\[
[[\phi,\phi]^{A^{0,2}(T^{1,0})},\bar{\psi}]^{\wedge(4)}=[[\phi,\phi]^{A^{1,1}(T^{1,0})},\bar{\psi}]^{\wedge(3)}=[[\phi,\phi]^{A^{0,2}(T^{0,1})},\bar{\psi}]^{\wedge(3)}=0,
\]
which implies that $[i_{[\phi,\phi]},i_{\bar{\psi}}]^{(4)}=0$.

For (5), first by (\ref{li commutator}) we have
$[\mathcal{L}_{\phi},i_{\overline{\psi}}]=i_{[\phi,\overline{\psi}]}-\mathcal{L}_{i_{\overline{\psi}}\phi}$. And so
\begin{equation*}
[\mathcal{L}_{\phi},i_{\overline{\psi}}]^{(2)}=[i_{[\phi,\overline{\psi}]}-\mathcal{L}_{i_{\overline{\psi}}\phi},i_{\overline{\psi}}]
=i_{[[\phi,\overline{\psi}],\overline{\psi}]^{\wedge} }-i_{[i_{\overline{\psi}}\phi,\overline{\psi}]},
\end{equation*}
where we have used the fact that $\mathcal{L}_{i_{\overline{\psi}}i_{\overline{\psi}}\phi}=0$ since $i_{\overline{\psi}}i_{\overline{\psi}}\phi=0$. Similar to \eqref{notice}, by~\cite[pp.\,70,\,Thm\,8.7\,,(6)]{KMS93} we have
\[
[\phi,\bar{\psi}]\in A^{1,1}(M,T^{1,0})\oplus A^{1,1}(M,T^{0,1})\oplus A^{0,2}(M,T^{0,1})\oplus A^{2,0}(M,T^{1,0}),
\]
and
\begin{align*}
&[[\phi,\bar{\psi}]^{A^{1,1}(T^{1,0})},\bar{\psi}]^{\wedge(3)}=[[\phi,\bar{\psi}]^{A^{1,1}(T^{0,1})},\bar{\psi}]^{\wedge(2)}=0,\\
&[[\phi,\bar{\psi}]^{A^{0,2}(T^{0,1})},\bar{\psi}]^{\wedge(3)}=[[\phi,\bar{\psi}]^{A^{2,0}(T^{1,0})},\bar{\psi}]^{\wedge(2)}=0,
\end{align*}
which implies that $[[\phi,\bar{\psi}],\bar{\psi}]^{\wedge(3)}=0$. Similarly, $[[i_{\bar{\psi}}\phi,\bar{\psi}],\bar{\psi}]^{\wedge(3)}=0$. Hence
\begin{align*}
&[\mathcal{L}_{\phi},i_{\overline{\psi}}]^{(3)}=i_{\big[[[\phi,\bar{\psi}] \bar{\psi}]^{\wedge} -[i_{\bar{\psi}}\phi,\bar{\psi}], \bar{\psi}\big]^{\wedge}}
=i_{[[\phi,\bar{\psi}] \bar{\psi}]^{\wedge(2)} -[[i_{\bar{\psi}}\phi,\bar{\psi}], \bar{\psi}]^{\wedge}}\\
&[\mathcal{L}_{\phi},i_{\overline{\psi}}]^{(4)}
=i_{[[\phi,\bar{\psi}] \bar{\psi}]^{\wedge(3)} -[[i_{\bar{\psi}}\phi,\bar{\psi}], \bar{\psi}]^{\wedge(2)}}= -i_{[[i_{\bar{\psi}}\phi,\bar{\psi}], \bar{\psi}]^{\wedge(2)}}\\
&[\mathcal{L}_{\phi},i_{\overline{\psi}}]^{(5)}=-i_{[[i_{\bar{\psi}}\phi,\bar{\psi}], \bar{\psi}]^{\wedge(3)}}=0.
\end{align*}

For (6), first we have
\begin{align*}
&e^{-i_{\overline{\psi}}}e^{-i_{\phi}}\nabla e^{i_{\phi}}e^{i_{\overline{\psi}}}\\
=&e^{-i_{\overline{\psi}}}(\nabla - \mathcal{L}_{\phi} -\frac{1}{2}i_{[\phi,\phi]}- \frac{1}{3!}i_{[[\phi,\phi],\phi]^{\wedge}})e^{i_{\overline{\psi}}}\\
=&\nabla-\mathcal{L}_{\overline{\psi}}-\frac{1}{2}i_{[\overline{\psi},\overline{\psi}]} - \frac{1}{3!}i_{[[\bar{\psi},\bar{\psi}],\bar{\psi}]^{\wedge}} -e^{-i_{\overline{\psi}}} \mathcal{L}_{\phi} e^{i_{\bar{\psi}}}-\frac{1}{2}e^{-i_{\overline{\psi}}}i_{[\phi,\phi]}e^{i_{\bar{\psi}}}- \frac{1}{3!}e^{-i_{\bar{\psi}}}i_{[[\phi,\phi],\phi]^{\wedge}}e^{i_{\bar{\psi}}},
\end{align*}
then (6) follows from (4),(5) and the fact that $\Big[ [[\phi,\phi],\phi]^{\wedge}, \bar{\psi} \Big]^{\wedge(4)}=0$.
\end{proof}

\begin{remark}It is not hard to check that the left hand side of these identities are all graded derivations. Hence we can also prove these formulas by using (\ref{D=lK+iL}). In a subsequent paper, we will carry out this approach to prove several extension formulas and study its applications in deformation of complex structures.
\end{remark}
\begin{remark}In the case $E$ is trivial, it's easy to see that $\mathcal{L}_{\phi}^{0,1}=[i_\phi,\overline{\partial}]$ is an algebraic derivation, that is, $\mathcal{L}_{\phi}^{0,1}f=0, \forall f\in A^0(M)$. Hence $\mathcal{L}_{\phi}^{0,1}=i_K$ for some unique $K\in A^{0,2}(M,T^{1,0})$. Let $\{e_i\}$ be a local frame of $T^{1,0}$ and $\{\xi^i\}\subset A^{1,0}(M)$ its dual frame, since $i_K\bar{\xi^i}=0$ and $i_K\xi^i=i_\phi\overline{\partial}\xi^i-\overline{\partial}i_\phi\xi^i$, we see that $K=K^i\otimes e_i$, where $K^i=i_\phi\overline{\partial}\xi^i-\overline{\partial}i_\phi\xi^i$. If furthermore $J$ is integrable, i.e. $M$ is a complex manifold, then we may set $e_i=\frac{\partial}{\partial z^i}$ and $\xi^i=dz^i$ to be the coordinate frame, it follows that $K=-\overline{\partial}\phi$ and $\mathcal{L}_{\phi}^{0,1}=-i_{\overline{\partial}\phi}$.
\end{remark}
\begin{remark}\label{rk-lie-alg}
In the case $E$ is trivial, it's easy to see that $\mathcal{L}_{\phi}^{0,1}=[i_\phi,\overline{\partial}]$ is an algebraic derivation, that is, $\mathcal{L}_{\phi}^{0,1}f=0, \forall f\in A^0(M)$. Hence $\mathcal{L}_{\phi}^{0,1}=i_K$ for some unique $K\in A^{0,2}(M,T^{1,0})$. Let $\{e_i\}$ be a local frame of $T^{1,0}$ and $\{\xi^i\}\subset A^{1,0}(M)$ its dual frame, since $i_K\bar{\xi^i}=0$ and $i_K\xi^i=i_\phi\overline{\partial}\xi^i-\overline{\partial}i_\phi\xi^i$, we see that $K=K^i\otimes e_i$, where $K^i=i_\phi\overline{\partial}\xi^i-\overline{\partial}i_\phi\xi^i$. If furthermore $J$ is integrable, i.e. $M$ is a complex manifold, then we may set $e_i=\frac{\partial}{\partial z^i}$ and $\xi^i=dz^i$ to be the coordinate frame, it follows that $K=-\overline{\partial}\phi$ and $\mathcal{L}_{\phi}^{0,1}=-i_{\overline{\partial}\phi}$.
\end{remark}

Since $(i_{\phi})^n=0$, by (4), we see that $(i_{\phi+[\phi,\overline{\psi}]^{\wedge}+[\phi,\overline{\psi}]^{\wedge(2)}})^n=0$. We failed to give a direct proof of this simple fact. Indeed, we have the following more general result
\begin{prop}
Let $R$ be an unitary associative algebra over $\mathbb{Q}$ (not necessarily commutative), $x,y\in R$ be nilpotent elements. Assume $x^N=0$ for some $N>1$ and $x\in R$ is $k$-commutable with $y$, i.e. $[x,y]^{(k)}=0$, then $(\sum_{i=0}^{i=k-1}\frac{1}{i!}[x,y]^{(i)})^N=0$ and
\begin{equation*}
e^{-y}e^xe^{y}=e^{\sum_{i=0}^{i=k-1}\frac{1}{i!}[x,y]^{(i)}}
\end{equation*}
\end{prop}

\begin{proof}This follows immediately from Lemma \ref{commutator lemma}.
\end{proof}

\vskip 1\baselineskip \textbf{Acknowledgements.} This note is motivated by the work \cite{LRY15},\cite{RZ17} of Liu-Rao-Yang and Rao-Zhao on extension formulas. I am grateful to Professor Kefeng Liu for his constant support and encouragement. I would like to thank Kang Wei and Hai-Sheng Liu for useful discussions. I would also like to thank the anonymous referees for their valuable suggestions.

\bibliographystyle{alpha}

\begin{thebibliography}{KMS93}

\bibitem[Dem12]{Dem12}
J.~P. Demailly.
\newblock {\em Complex analytic and differential geometry}.
\newblock 2012.
\newblock available
  at~~{https://www-fourier.ujf-grenoble.fr/demailly/manuscripts/agbook.pdf}.

\bibitem[FL19]{FL19}
J.-X. Fu and H.-S. Liu.
\newblock Extension formulae on almost complex manifolds.
\newblock arXiv:1903.09821, 2019.

\bibitem[Got05]{Got05}
R.~Goto.
\newblock On deformations of generalized {Calabi-Yau}, {hyperK\"ahler, G2 and
  Spin(7)} structures {I}.
\newblock arXiv:math/0512211v1, 2005.

\bibitem[KMS93]{KMS93}
I.~Kol{\'a}$\check{r}$, P.~W. Michor, and J.~Slov{\'a}k.
\newblock {\em Natural operations in differential geometry}.
\newblock Springer-Verlag, Berlin, Heidelberg, New York, 1993.

\bibitem[LR11]{LR11}
K.~Liu and S.~Rao.
\newblock Remarks on the {C}artan formula and its applications.
\newblock {\em Asian Journal of Mathematics}, 16(1):157--169, 2011.

\bibitem[LRY15]{LRY15}
K.~Liu, S.~Rao, and X.~Yang.
\newblock Quasi-isometry and deformations of {C}alabi-{Y}au manifolds.
\newblock {\em Inventiones mathematicae}, 199(2):423--453, 2015.

\bibitem[Man05]{Ma05}
M.~Manetti.
\newblock Lectures on deformations of complex manifolds.
\newblock arXiv:math/0507286v1 [math.AG], 2005.

\bibitem[Mic86]{Mic86}
P.~W. Michor.
\newblock Remarks on the {F}r{\"o}licher-{N}ijenhuis bracket.
\newblock In {\em Proceedings of the Conference on Differential Geometry and
  its Applications, Brno}, volume 1086, pages 197--220, 1986.

\bibitem[MK71]{MK71}
J.~Morrow and K.~Kodaira.
\newblock {\em Complex manifolds}, volume 355.
\newblock American Mathematical Soc., 1971.

\bibitem[RZ17]{RZ17}
S.~Rao and Q.~Zhao.
\newblock Several special complex structures and their deformation properties.
\newblock The Journal of Geometric Analysis, 2017.
\newblock {https://doi.org/10.1007/s12220-017-9944-7}.

\bibitem[Xia19]{Xia19dp}
W.~Xia.
\newblock Derivations on almost complex manifolds.
\newblock {\em Proceedings of the American Mathematical Society}, 147:559--566,
  2019.

\end{thebibliography}

\end{document}